\documentclass[11pt]{article}
\usepackage{amsmath,amssymb,theorem}
\usepackage{graphicx,subfigure,enumerate}
\usepackage[ruled, lined, linesnumbered]{algorithm2e}
\usepackage{psfrag,placeins,wrapfig,booktabs,chngcntr}
\counterwithin{figure}{section}
\numberwithin{figure}{section}
\counterwithin{table}{section}
\usepackage{hyperref}
\hypersetup{
  colorlinks,
  citecolor=black,
  filecolor=black,
  linkcolor=black,
  urlcolor=black
}

\newtheorem{thm}{Theorem}[section]
\newtheorem{conj}[thm]{Conjecture}
\newtheorem{lem}[thm]{Lemma}
\theorembodyfont{\rmfamily}
\def\pf{\bigskip\noindent {\bf Proof.}~~}

\def\dfn#1{{\sl #1}}
\def\es{\emptyset}
\def\less{\setminus}

\def\pf{\bigskip\noindent {\bf{Proof.}}~~}

\newcounter{counter}

  \allowdisplaybreaks[4]

\begin{document}
\title{On the size  of  $(K_t,\mathcal{T}_k)$-co-critical  graphs}
\author{Zi-Xia Song\thanks{Supported by the National   Science  Foundation under Grant No. DMS-1854903. }\,  and Jingmei Zhang\thanks{E-mail addresses: Zixia.Song@ucf.edu (Z-X. Song); jmzhang@Knights.ucf.edu (J. Zhang)}  \\
Department of Mathematics\\
University of Central Florida\\
Orlando, FL 32816, USA\\
}
 \date{ }
\maketitle

\begin{abstract}
Given an integer $r\ge1$ and graphs $G, H_1, \ldots,  H_r$, we write \emph{$G \rightarrow ({H}_1, \ldots,  {H}_r)$} if  every $r$-coloring of the edges of $G$ contains  a monochromatic copy of $H_i$ in color $i$ for some $i\in\{1, \ldots, r\}$. A non-complete graph $G$ is \emph{$(H_1, \ldots,  H_r)$-co-critical} if $G  \nrightarrow ({H}_1,  \ldots, {H}_r)$,   but  $G+e\rightarrow ({H}_1,  \ldots, {H}_r)$ for every edge $e$ in $\overline{G}$.  In this paper,  motivated  by Hanson and Toft's conjecture   [Edge-colored saturated graphs, J  Graph Theory  11(1987), 191--196],   we study the minimum number of edges over all  $(K_t, \mathcal{T}_k)$-co-critical graphs on $n$ vertices, where $\mathcal{T}_k$ denotes the family of all trees on $k$ vertices.      Following  Day [Saturated graphs of prescribed minimum degree,  Combin.  Probab.   Comput.  26 (2017),   201--207], we  apply   graph bootstrap percolation on a not necessarily $K_t$-saturated graph  to   prove   that for all   $t\ge4 $ and      $k\ge \max\{6, t\}$,    there exists a constant $c(t, k)$ such that, for all      $n \ge (t-1)(k-1)+1$, if $G$  is  a $(K_t, \mathcal{T}_k)$-co-critical graph  on $n$ vertices, then    $$   e(G)\ge   \left(\frac{4t-9}{2}+\frac{1}{2}\left\lceil \frac{k}{2} \right\rceil\right)n-c(t, k).$$ 
Furthermore, this linear bound is asymptotically best possible  when $t\in\{4,5\}$ and   $k\ge6$. The method we develop in this paper may shed some light on attacking Hanson and Toft's conjecture.  
 
\end{abstract}


\baselineskip 16pt

\section{Introduction}
All graphs considered in this paper are finite, and without loops or multiple edges. For a graph $G$, we will use $V(G)$ to denote the vertex set, $E(G)$ the edge set, $|G|$ the number of vertices, $e(G)$ the number of edges,  $N_G(x)$   the neighborhood of vertex $x$ in $G$, $\delta(G)$ the minimum degree, $\Delta(G)$ the maximum degree,  and $\overline{G}$ the complement of $G$.
If  $A, B \subseteq V(G)$ are disjoint, we say that $A$ is \dfn{complete to} $B$ if  every vertex in $ A$ is adjacent to every vertex in   $B$; and $A$ is  \dfn{anti-complete to} $B$ if   no vertex in $A$ is adjacent to a vertex in $B$.
The subgraph of $G$ induced by $A$, denoted $G[A]$, is the graph with vertex set $A$ and edge set $\{xy \in E(G): x, y \in A\}$. We denote by $B \less A$ the set $B - A$, $e_G(A, B)$ the number of edges between $A$ and $B$ in $G$, and $G \less A$ the subgraph of $G$ induced on $V(G) \less A$, respectively.
If $A = \{a\}$, we simply write $B \less a$, $e_G(a, B)$, and $G \less a$, respectively.  For  any  edge  $e$ in $ \overline{G} $, we use $G+e$ to denote the graph obtained from $G$ by adding the new edge $e$. 
The {\dfn{join}} $G+H$ (resp.~{\dfn{union}} $G\cup H$) of two 
vertex disjoint graphs
$G$ and $H$ is the graph having vertex set $V(G)\cup V(H)$  and edge set $E(G)
\cup E(H)\cup \{xy:   x\in V(G),  y\in V(H)\}$ (resp. $E(G)\cup E(H)$).
Given two isomorphic graphs $G$ and $H$, we may (with a slight but common abuse of notation) write $G = H$.   For an integer $t\ge1$ and a graph $H$, we define $tH$ to be the union of $t$ disjoint copies of $H$.  We use $K_n$, $K_{1,{n-1}}$, $C_n$, $P_n$ and $T_n$ to denote the complete graph,  star,  cycle,   path and a tree on  $n$ vertices, respectively.  For any positive integer $r$, we write  $[r]$ for the set $\{1,2, \ldots, r\}$. We use the convention   ``$A:=$'' to mean that $A$ is defined to be the right-hand side of the relation.  \medskip

Given an integer $r \ge 1$ and   graphs $G$, ${H}_1, \dots, {H}_r$, we write \dfn{$G \rightarrow ({H}_1, \dots, {H}_r)$} if every $r$-coloring of $E(G)$ contains a monochromatic  ${H}_i$ in color $i$ for some $i\in [r]$.
The classical \dfn{Ramsey number}  $R({H}_1, \dots, {H}_r)$  is the minimum positive integer $n$ such that $K_n \rightarrow ({H}_1, \dots, {H}_r)$.  Following Ne\v{s}et\v{r}il~\cite{Nesetril1986}, and Galluccio, Siminovits and Simonyi~\cite{Galluccio1992}, we say that  a non-complete graph $G$   is \emph{$(H_1, \ldots,  H_r)$-co-critical} if $G  \nrightarrow ({H}_1,  \ldots, {H}_r)$,   but  $G+e\rightarrow ({H}_1,  \ldots, {H}_r)$ for every edge $e$ in $\overline{G}$.  Clearly,  $K_6^-$ is $(K_3, K_3)$-co-critical, where   $K_6^-$ denotes  the graph obtained from $K_6$ by deleting exactly one edge. 
It is worth noting that every   $(H_1, \ldots,  H_r)$-co-critical graph has at least $R({H}_1, \dots, {H}_r)$ many vertices.  \medskip

\noindent{\bf Remark}.  Following Galluccio, Siminovits and Simonyi~\cite{Galluccio1992}, we excluded the complete graphs in the definition of $(H_1, \ldots,  H_r)$-co-critical graphs, else every complete graph on fewer than $R({H}_1, \dots, {H}_r)$ vertices is  $(H_1, \ldots,  H_r)$-co-critical.  \medskip

The notation of co-critical graphs was initiated by Ne\v{s}et\v{r}il~\cite{Nesetril1986} in 1986 when he asked the following question regarding   $(K_3, K_3)$-co-critical graphs: 
\begin{quote}
 Are there   infinitely  many \dfn{minimal} co-critical graphs, i.e.,  co-critical graphs which lose this property when any vertex is deleted? Is $K_6^-$ the only one? 
\end{quote}
This was answered in the positive by Galluccio, Siminovits and Simonyi~\cite{Galluccio1992}.  They constructed infinite many minimal $(K_3, K_3)$-co-critical graphs without containing $K_5$ as a subgraph.    Szab\'o~\cite{Szabo1996} then constructed   infinitely many  nearly regular $(K_3, K_3)$-co-critical graphs with low maximum degree.  It remains unknown whether there are infinitely many \dfn{strongly} minimal co-critical graphs, where an $(H_1, \ldots,  H_r)$-co-critical graph is \dfn{strongly minimal co-critical} if it contains no proper subgraph which is also $(H_1, \ldots,  H_r)$-co-critical. This is one of the most   intriguing open problems proposed by   Galluccio, Siminovits and Simonyi in \cite{Galluccio1992}.   One interesting observation made in   \cite{Galluccio1992} is that  if $G$ is $(H_1, \ldots,  H_r)$-co-critical, then $\chi(G)\ge R({H}_1, \dots, {H}_r)-1$.  They also made some observations on the minimum  degree of $(K_3, K_3)$-co-critical graphs and maximum number of possible edges of $(H_1, \ldots,  H_r)$-co-critical graphs.  
 \medskip

We want to point out here that Hanson and Toft~\cite{Hanson1987} in 1987 also studied  the minimum and maximum number  of edges over all $(H_1, \ldots,  H_r)$-co-critical graphs on $n$ vertices when $H_1, \ldots,  H_r$ are complete graphs, under the name of \dfn{strongly $(|H_1|, \ldots, |H_r|)$-saturated} graphs. Recently, this topic has  been studied under the name of \dfn{$\mathcal{R}_{\min}(H_1, \dots, H_r)$-saturated} graphs \cite{Chen2011, Ferrara2014,  RolekSong}. A graph $G$ is \dfn{$({H}_1, \dots, {H}_r)$-Ramsey-minimal} if $G \rightarrow ({H}_1, \dots, {H}_r)$, but for any proper subgraph $G'$ of $G$, $G' \nrightarrow ({H}_1, \dots, {H}_r)$.
We define  $\mathcal{R}_{\min}({H}_1, \dots, {H}_r)$ to be  the family of all $({H}_1, \dots, {H}_r)$-Ramsey-minimal graphs.  A graph $G$ is \dfn{$\mathcal{R}_{\min}(H_1, \dots, H_r)$-saturated} if  no element  of $\mathcal{R}_{\min}(H_1, \dots, H_r)$ is a subgraph of $G$, but for any edge $e$ in $\overline{G}$,  some element of $\mathcal{R}_{\min}(H_1, \dots, H_r)$ is  a subgraph of $G + e$. It can be easily checked that a non-complete graph is $(H_1, \ldots,  H_r)$-co-critical if and only if it is 
$\mathcal{R}_{\min}(H_1, \dots, H_r)$-saturated.  From now on, we shall  use  the notion of   $(H_1, \ldots,  H_r)$-co-critical  other than       $\mathcal{R}_{\min}(H_1, \dots, H_r)$-saturated, as the former  is  much simpler and   straightforward.  \medskip

Let   $R = R(K_{t_1}, \dots, K_{t_r})$ be the classical Ramsey number for $K_{t_1}, \dots, K_{t_r}$. Hanson and Toft~\cite{Hanson1987} proved that every $(K_{t_1}, \dots, K_{t_r})$-co-critical on $n$ vertices has at most $e(T_{R-1, n})$ edges and this bound is best possible, where $T_{R-1, n}$ is the Tur\'an graph on $n$ vertices without $K_{R}$. They also observed that for all $n\ge R$, the  graph $K_{R-2}+\overline{K}_{n-R+2}$ is $(K_{t_1}, \dots, K_{t_r})$-co-critical.  They further  made the following conjecture that no   $(K_{t_1}, \dots, K_{t_r})$-co-critical graph on $n$ vertices can have fewer than $e(K_{R-2}+\overline{K}_{n-R+2})$ edges.

\begin{conj}[Hanson and Toft~\cite{Hanson1987}]\label{HTC}  
Let $G$ be a  $(K_{t_1}, \dots, K_{t_r})$-co-critical graph   on $n$ vertices. Then   
$$e(G)\ge  (R- 2)(n - R + 2) + \binom{R - 2}{2}.$$
This bound is best possible for every $n$. 
\end{conj}
\medskip

Conjecture~\ref{HTC} remains wide open, except that the first nontrivial case,  $(K_3,  K_3)$-co-critical graphs,   has        been settled in \cite{Chen2011} for $n \ge 56$.  Structural properties of $(K_3,  K_4)$-co-critical graphs are given in \cite{K3K4}.    Motivated by Conjecture~\ref{HTC},  we study  the following problem.  Let $\mathcal{T}_k$ denote  the family of all trees on $k$ vertices.
For all $t, k \ge 3$, we write $G\rightarrow (K_t, \mathcal{T}_k)$  if  for every $2$-coloring   $\tau: E(G) \to \{\text{red, blue} \}$, $G$ has either a red $K_t$ or a blue tree $T_k\in \mathcal{T}_k$. 
A non-complete graph $G$ is \dfn{$(K_t, \mathcal{T}_k)$-co-critical} if  
$G\nrightarrow (K_t, \mathcal{T}_k)$, but $G+e\rightarrow (K_t, \mathcal{T}_k)$ for all $e$ in $\overline{G}$.  The main purpose of this paper is to study the structural properties of      $ (K_t, \mathcal{T}_k)$-co-critical graphs on $n$ vertices in order to obtain the minimum size among all such graphs.  By a classic result of Chv\'atal~\cite{Chvatal}, $R(K_t, \mathcal{T}_k)=(t-1)(k-1)+1$.  Hence,   every $ (K_t, \mathcal{T}_k)$-co-critical graph has at least $ R(K_t, \mathcal{T}_k) =(t-1)(k-1)+1$ many vertices.  Following the observation made in both \cite{Galluccio1992} and  \cite{Hanson1987},   every $ (K_t, \mathcal{T}_k)$-co-critical graph  on $n$ vertices has at most $e(T_{R(K_t, \mathcal{T}_k)-1, n})$ edges.  
 We focus on studying the minimum number of possible edges over all $ (K_t, \mathcal{T}_k)$-co-critical graphs on $n$ vertices.  Very recently,  Rolek and the first author~\cite{RolekSong} proved the following.

\begin{thm}[Rolek and Song~\cite{RolekSong}]\label{K3Tk}  Let $n, k\in \mathbb{N}$.  
\begin{enumerate}[(i)]
 \item  Every  $(K_3, \mathcal{T}_4)$-co-critical graph on $n\ge 18$ vertices has  at least $\left\lfloor 5n/2\right\rfloor$ edges. This bound is sharp for every $n\ge18$.  
 \item   For all  $k \ge 5$,  if $G$ is $(K_3, \mathcal{T}_k)$-co-critical on $n\ge  2k + (\lceil k/2 \rceil +1) \lceil k/2 \rceil -2$ vertices, then $$ e(G) \ge \left(\frac{3}{2}+\frac{1}{2}\left\lceil \frac{k}{2} \right\rceil\right)n-c(k),$$ 
 where $c(k)=\left(\frac{1}{2} \left\lceil \frac{k}{2} \right\rceil + \frac{3}{2} \right) k -2$. This bound is asymptotically best possible.
 \end{enumerate}
\end{thm}

To state our results, we need to introduce more notation.  Given a family $\mathcal{F}$, a graph   is \dfn{$\mathcal{F}$-free} if it does not contain any graph $F\in \mathcal{F}$ as a subgraph. We simply say a graph is \dfn{$F$-free} when $\mathcal{F}=\{F\}$.  Erd\H os, Hajnal   and Moon~\cite{Erdos1964} in 1964 initiated the study of the minimum number of edges over all $K_t$-saturated graphs on $n$ vertices (see the dynamic survey \cite{Faudree2011} on  the extensive studies on $K_t$-saturated graphs).  Theorem \ref{Day} below is a result of Day \cite{Day2017} on $K_t$-saturated graphs with prescribed minimum degree. It confirms a conjecture of Bollob\'as \cite{Bollobas} when $t = 3$.  It is worth noting that Day applied    the  $r$-neighbour    bootstrap percolation on  a $K_t$-saturated graph  to prove Theorem~\ref{Day}, where  graph bootstrap percolation was introduced in \cite{bootstrap}.  
Theorem~\ref{Hajnal} is a result of Hajnal   \cite{Hajnal}   
on $K_t$-saturated graphs. It seems hard to improve  $2(t-2)$ in Theorem~\ref{Hajnal}. 
\begin{thm}[Day \cite{Day2017}]\label{Day} 
Let $q \in\mathbb{N}$. There exists a constant $c = c(q)$ such that, for  all $3\le t \in \mathbb{N}$ and all $n\in \mathbb{N}$, if $G$ is a $K_t$-saturated graph on $n$ vertices with  $\delta(G) \ge q$, then $e(G) \ge qn - c$.
\end{thm}

\begin{thm}[Hajnal  \cite{Hajnal}]\label{Hajnal}
Let $t, n \in\mathbb{N}$. Let $G$ be a $K_t$-saturated graph on $n$ vertices.  Then either $\Delta(G) = n-1$ or $\delta(G) \ge 2(t-2)$. 
\end{thm}

For a $ (K_t, \mathcal{T}_k)$-co-critical graph $G$,  let   $\tau : E(G) \to \{\text{red, blue} \}$  be a  $2$-coloring of $E(G)$ and let  $E_r$ and $E_b$ be the  color classes of the coloring $\tau$.  We use $G_{r}$ and $G_{b}$ to denote the spanning subgraphs of $G$  with edge sets  $E_r$ and $E_b$, respectively.  
We define  $\tau$ to be a  \dfn{critical-coloring} of $G$ if $G$ has neither  a red $K_t$ nor a blue $T_k\in  \mathcal{T}_k$ under $\tau$, that is,  if  $G_r$ is $K_t$-free and $G_b$ is $\mathcal{T}_k$-free.  For every $v\in V(G)$, we use $d_r(v)$ and $N_r(v)$ to denote the degree and neighborhood of $v$ in $G_r$, respectively. Similarly, we define $d_b(v)$ and $N_b(v)$ to be the degree and neighborhood of $v$ in $G_b$, respectively. One can see that if  $G$ is $ (K_t, \mathcal{T}_k)$-co-critical, then $G$ admits at least one critical-coloring but $G+e$ admits no critical-coloring for every edge $e$ in $\overline{G}$.\medskip 

In this paper, we  first  establish  a number of important structural properties of  $(K_t, \mathcal{T}_k)$-co-critical graphs in the hope that the method we develop here may shed some light on attacking Conjecture~\ref{HTC}.     Theorem \ref{e(H)}(\ref{eqH}) below  is crucial in the proof of Theorem~\ref{lower}.   
 Following Day  \cite{Day2017},  we   apply the $q$-neighbour bootstrap percolation on a    not necessarily $K_t$-saturated graph, to prove Theorem \ref{e(H)}(\ref{eqH}), but with more involved   rules.   
     \medskip

\begin{thm}\label{e(H)}
For all $t, k\in\mathbb{N}$ with  $t \ge 3$ and $k \ge 3$, let $G$ be a  $(K_t, \mathcal{T}_k)$-co-critical graph on  $n$   vertices.  Among all critical-colorings of $G$, let  $\tau : E(G) \rightarrow \{$red, blue$\}$ be  a critical-coloring of $G$ with $|E_r|$  maximum. Let $D_1, \ldots, D_p$ be all components of $G_b$. Let $H:=G \less (\bigcup_{i \in [p]} E(G[V(D_i)]))$.  Then the following hold.
\begin{enumerate}[\rm(a)]
\item\label{Gr} $\Delta(G_r) \le n-2$ and $\delta(G_r) \ge 2(t-2)$.

\item\label{NHu} For all  $i, j \in [p]$ with $i\ne j$, if there exist  $u \in V(D_i)$ and   $v \in V(D_j)$  such that   $uv \notin E(H)$, then   $H[N_H(u) \cap N_H(v)]$ contains  a $K_{t-2}$  subgraph.

\item\label{Di}  For every   $uv \in E(H)$, if  $v$ is contained in all  $K_{t-2}$ subgraphs of    $H[N_H(u)]$ and $\{v\}=V(D_j)$ for some $j\in[p]$, then   $|D_i| =k-1$ for all $D_i$ with $u\notin D_i$ and $ D_i \less N_H(u) \ne  \emptyset$, where  $i \in [p]$.  

\item\label{nocom} If $\delta(H) \le 2t-5$ and $k \ge t$, then for any vertex   $u \in V(H)$ with $d_H(u)=\delta(H)$,  no edge  of $H[N_H(u)]$ is contained in all   $K_{t-2}$ subgraphs of    $H[N_H(u)]$.

\item\label{dH} $k\ge 2t-1-\delta(H)$ and $\delta(H) \ge t-1$.

\item\label{eD}    $\displaystyle \sum_{i=1}^p e(G[V(D_i)]) > \left( \frac{1}{2}\left\lceil \frac{k}{2} \right\rceil -\frac{1}{2}\right)(n-(t-1)(\lceil k/2 \rceil-1))$. 

\item\label{conH}  $H$ is connected. 

\item\label{eqH}  For every $q\in\mathbb N$ with $q \ge t-1$,  there exists a constant $ c(q,k)$ such that,    if $\delta(H) \ge q $, then $e(H) \ge qn-c(q, k)$.

\end{enumerate}
\end{thm}

 We prove Theorem~\ref{e(H)} in Section~\ref{property}. We then  apply Theorem~\ref{e(H)} to  study  the size of $(K_t, \mathcal{T}_k)$-co-critical graphs.  We prove   Theorem~\ref{lower} in Section~\ref{Lower}.   
 
 \begin{thm}\label{lower}
 Let  $ t, k\in \mathbb{N}$ with $t \ge 4$ and  $k\ge\max\{6, t\}$.  There exists a  constant  $\ell(t, k)$  such that, for  all $n\in \mathbb{N}$ with    $n \ge (t-1)(k-1)+1$,    if $G$  is  a $(K_t, \mathcal{T}_k)$-co-critical graph  on $n$ vertices, then    
 $$   e(G)\ge   \left(\frac{4t-9}{2}+\frac{1}{2}\left\lceil \frac{k}{2} \right\rceil\right)n-\ell(t, k).$$ 
\end{thm}

Finally we prove that   the linear bound  given in Theorem~\ref{lower} is asymptotically best possible    when $t\in\{4,5\}$ and $k\ge6$.    Proof of Theorem~\ref{upper} is given in  Section~\ref{Upper}. 

\begin{thm}\label{upper}
For each $t \in \{4, 5\}$,  all $k \ge 3$ and  $n \ge (2t-3)(k-1)+{\lceil k/2 \rceil}{\lceil k/2 \rceil}-1$, there exists a   $(K_t, \mathcal{T}_k)$-co-critical graph  $G$ on $n$ vertices such that 
$$e(G) \le \left(\frac{4t-9}{2}+\frac{1}{2} \left\lceil \frac{k}{2} \right\rceil \right)n+C(t, k)$$
where $C(t, k)=\frac{1}{2}(t^2+t-5)k^2-(2t^2+2t-11)k  -\frac{(t-2)(t-19)}{2} 
- \frac{1}{2} \left\lceil \frac{k}{2} \right\rceil \left((2t-3)(k-1) - \left\lceil \frac{k}{2} \right\rceil  \right)$.
\end{thm}

 With the support of Theorem~\ref{K3Tk} and  Theorem~\ref{upper}, we believe 
 that the linear bound given  in Theorem~\ref{lower} is  
  asymptotically best possible for all $t\ge3$ and $k\ge3$.

\section{Structural properties of  $(K_t, \mathcal{T}_k)$-co-critical graphs}\label{property}

We first prove  the following lemma.

\begin{lem}\label{structure}
For all $t, k\in\mathbb{N}$ with  $t \ge 3$ and $k \ge 3$, let $G$ be a  $(K_t, \mathcal{T}_k)$-co-critical graph on $n $ vertices.  Let $\tau : E(G) \rightarrow \{$red, blue$\}$ be a critical-coloring of $G$. Then the following hold. 
\begin{enumerate}[(a)]
\item For every component  $D$   of $G_b$,  $|D|\le k-1$ and $G[V(D)]=K_{|D|}$.

\item If $D_1, \cdots, D_q$ are the components of $G_b$ with $|D_i| <k/2$ for all $i \in [q]$,  then $V(D_1), \cdots, V(D_q)$ are complete to each other in $G_r$, and so $q \le t-1$.
\end{enumerate}
\end{lem}

\pf 
To prove (a), let $D$ be a  component  of $G_b$.  Since $G_b$ is $\mathcal{T}_k$-free, we see that $|D|\le k-1$. Suppose next that  $G[V(D)]\neq K_{|D|}$. Let $u,v \in V(D)$ be such that $uv \not\in E(G)$. We obtain a critical-coloring  of $G+uv$ from  $\tau$ by coloring the edge $uv$ blue, a contradiction.\medskip

To prove (b), let $D_1, \cdots, D_q$ be the components of $G_b$ with $|D_i| <k/2$ for all $i \in [q]$. Since $G$ is $(K_t, \mathcal{T}_k)$-co-critical, we see that $G+e$ admits no critical-coloring for every edge   $e$ in $\overline{G}$.   Let  $i, j \in [q]$ with $i \neq j$. We next  show that  $V(D_i)$ is complete to $V(D_j)$ in $G_r$. Suppose that there exist vertices $u \in V(D_i)$ and $v \in V(D_j)$ such that $uv \not\in E_r$. Then $uv \not\in E(G)$ and so we obtain a critical-coloring  of $G+uv$ from  $\tau$ by coloring the edge $uv$ blue, a contradiction. Thus $V(D_i)$ is complete to $V(D_j)$ in $G_r$ for all $i,j \in [q]$ with $i \neq j$. Since $\tau$ is a critical-coloring,  it follows that $G_r$ is $K_t$-free and so  $q \le t-1$. 
\hfill\vrule height3pt width6pt depth2pt \bigskip

We are now ready to prove Theorem~\ref{e(H)}. \medskip

\noindent {\bf Proof of Theorem~\ref{e(H)}:} Let $G$, $\tau$, $D_1, \ldots, D_p$ and $H$ be given as in the statement.  Then $n\ge (t-1)(k-1)+1$.  By Lemma~\ref{structure}(a), $|D_i|\le k-1$ for all $i\in [p]$.  Hence,  $G_b$ has at least $t$ components because $|G_b|=n \ge (t-1)(k-1)+1$. We first  prove Theorem~\ref{e(H)}(\ref{Gr}).  
  By the choice of $\tau$, $G_r$ is $K_t$-free but $G_r + e$ contains a copy of  $K_t$ for every  $e \in E(\overline{G_r})$.  Hence $G_r$ is $K_t$-saturated.  
Suppose there exists a vertex  $x \in V(G)$ such that  $d_r(x) = n-1$.   Note that $G_r\less x$ is $K_{t-1}$-free because $G_r$ is $K_t$-free. Since $G\ne K_n$,  there must exist  $u, w \in N_r(x)$ such that  $uw \not\in E(G)$.  By Lemma~\ref{structure}(a), $u,w$ belong to different components of $G_b$.   
 But then we obtain a critical-coloring  of $G+uw$ from  $\tau$ by first coloring the edge $uw$ red, and then recoloring $xu$ blue and all   edges incident with $u$ in $G_b$ red, a contradiction. This proves that $\Delta(G_r) \le n-2$. Since $G_r$ is $K_t$-saturated, by Theorem \ref{Hajnal}, $\delta(G_r) \ge 2(t-2)$.\medskip

To prove Theorem~\ref{e(H)}(\ref{NHu}), let $u \in V(D_i)$ and   $v \in V(D_j)$  be such that   $uv \notin E(H)$, where $i\ne j$.   Suppose    $H[N_H(u) \cap N_H(v)]$ is  $K_{t-2}$-free. Since  $|D_\ell|\le k-1$ for all $\ell\in [p]$, we obtain a critical-coloring of $G+uv$ from $\tau$ by first coloring the edge $uv$ red, and then recoloring all red edges   in  $G[V(D_\ell)]$ blue  for all $\ell\in [p]$,  a contradiction. 
Therefore, $H[N_H(u) \cap N_H(v)]$ contains  a $K_{t-2}$ subgraph. This proves Theorem~\ref{e(H)}(\ref{NHu}).    \medskip

 To prove Theorem~\ref{e(H)}(\ref{Di}), let $uv \in E(H)$  be such that $v$ is contained in all  $K_{t-2}$ subgraphs of    $H[N_H(u)]$ and $\{v\}=V(D_j)$ for some $j\in[p]$. We may assume that $u \in V(D_p)$ and $\{v\}= V(D_{p-1})$.   Note that $H[N_H(u)]\less v$ is $K_{t-2}$-free.     
 Suppose there exists an $\ell\in[p-2]$ such that $ D_\ell\less N_H(u)\ne \emptyset$ but $|D_\ell| \le k-2$.   Let $w \in V(D_\ell)\less   N_H(u)$. Then $wv\in E_r$, else we   obtain a critical-coloring of $G+wv$ from  $\tau$ by coloring the edge $wv$  blue.  
 Since $H[N_H(u)]\less v$ is $K_{t-2}$-free, we  then  obtain a critical-coloring of $G+uw$ from  $\tau$ by coloring the edge $uw$ red, and then recoloring $wv$ blue and all red edges incident with $u$ in $G[V(D_p)]$ blue, a contradiction. This proves Theorem~\ref{e(H)}(\ref{Di}).  \medskip

To prove Theorem~\ref{e(H)}(\ref{nocom},\ref{dH}), let  $u\in V(H)$ with $d_H(u)=\delta(H) $. 
 We may assume that $u \in V(D_p)$.  By Theorem~\ref{e(H)}(\ref{NHu}), $d_H(u) \ge t-2$. Let $N_H(u):=\{u_1, \ldots,  u_{\delta(H)}\}$. By Theorem~\ref{e(H)}(\ref{NHu}) applied to $u$ and any vertex in  $V(H) \less (V(D_p) \cup N_H(u))$, we see that $ H[N_H(u)] $ must have  a $K_{t-2}$ subgraph.  We may  assume that $  H[\{u_1, \ldots, u_{t-2}\}]=K_{t-2} $. Then we may further assume that $u_i \in V(D_{p-i})$ for all $i \in [t-2]$.     Let $v \in V(H) \less (V(D_p) \cup N_H(u))$. \medskip

 To proceed to prove Theorem~\ref{e(H)}(\ref{nocom}), assume $d_H(u)\le 2t-5$ and $k \ge t$.     
  Suppose $H[N_H(u)]$ has an edge, say  $u_1u_2   $, that   is contained in all   $K_{t-2}$ subgraphs of    $H[N_H(u)]$.  Then both $H[N_H(u)]\less u_1$ and $H[N_H(u)]\less u_2$ are $K_{t-2}$-free.  By Theorem~\ref{e(H)}(\ref{NHu}) applied to $u$ and any vertex in $V(H) \less (V(D_p) \cup N_H(u))$, $V(H) \less (V(D_p) \cup N_H(u))$ must be complete to  $\{u_1, u_2\}$ in $H$.  Then $V(D_{p-1}) \cup V(D_{p-2}) \subseteq  N_H(u)\less \{u_3, \ldots, u_{t-2}\}$. Thus  $|V(D_{p-1}) \cup V(D_{p-2})| =\delta(H)-(t-4)\le t-1 \le k-1$,  because $\delta(H) \le 2t-5$ and $t \le k$.   Then we obtain a  critical-coloring of $G+uv$ from  $\tau$ by first coloring the edge $uv$ red, and then recoloring $u_1u_2$ blue and all red edges incident with $u$ in $G[V(D_p)]$ blue, a contradiction. This proves Theorem~\ref{e(H)}(\ref{nocom}).   \medskip
  
To proceed to prove Theorem~\ref{e(H)}(\ref{dH}), note that   $| N_r(u)\cap V(D_p) |  = |N_r(u)|-d_{H}(u)$. 
By Theorem~\ref{e(H)}(\ref{Gr}), $|N_r(u)|\ge 2t-4$.   Since $D_p$ is a component of $G_b$, we see that $N_b(u)\cap V(D_p)\ne \es$. It follows that $|V(D_p)|= |\{u\}|+|N_b(u)\cap V(D_p)|+| N_r(u)\cap V(D_p) | \ge 1+1+(2t-4)-d_{H}(u)=2t-2-d_{H}(u)$. 
By Lemma~\ref{structure}(a), $2t-2-d_{H}(u)\le | V(D_p) | \le k-1$, which yields  $k\ge 2t-1-d_{H}(u)$. 
Suppose next that $\delta(H) = t-2<2t-5$.  Then $k\ge t+1$. 
But then $  H[\{u_1, \ldots, u_{t-2}\}]$ is the only $K_{t-2}$ subgraph of $H[N_H(u)]$, contrary to Theorem~\ref{e(H)}(\ref{nocom}).  This proves Theorem~\ref{e(H)}(\ref{dH}).  \medskip

We next prove Theorem~\ref{e(H)}(\ref{eD}). 
 By Lemma \ref{structure}(a,b), $|D_i|\le k-1$,  $G[V(D_i)]=K_{|D_i|}$ for all $i \in [p]$, and  at most  $t-1$  of the  $D_i$'s have less than $k/2$ vertices. Let $r$ be the remainder of $n-(t-1)(\lceil k/2 \rceil-1)$ when divided by $\lceil k/2 \rceil$, and let $s \ge 0$ be an integer such that 
 \[n-(t-1)(\lceil k/2 \rceil-1) = s \lceil k/2 \rceil + r(\lceil k/2 \rceil +1).\]
  It is
straightforward to see that $\sum_{i=1}^p e(G[V(D_i)])$ is minimized when: $t-1$ of the components, say $D_1, \ldots, D_{t-1}$ are such that $|D_1|, \ldots, |D_{t-1}| < k/2$; $r$ of the components, say $D_t, \cdots, D_{r+t-1}$ are such that $|D_t|= \cdots =|D_{r+t-1}| = \lceil k/2 \rceil +1$; and $s$ of the components, say $D_{r+t}, \cdots, D_{r+s+t-1}$ are such that $|D_{r+t}|= \cdots =|D_{r+s+t-1}| = \lceil k/2 \rceil$.   Using the facts that  $s \lceil k/2 \rceil + r(\lceil k/2 \rceil +1)=n-(t-1)(\lceil k/2 \rceil-1)$ and $r\le \lceil k/2 \rceil -1$, it follows that 
\begin{align*}
\sum_{i=1}^p e(G[V(D_i)]) & > s{\lceil k/2 \rceil\choose 2} + r {\lceil k/2 \rceil +1 \choose 2}\\
&= \frac{s}{2} \left\lceil \frac{k}{2} \right\rceil  \left(\left\lceil \frac{k}{2} \right\rceil-1 \right)+\frac{r}{2} \left\lceil \frac{k}{2} \right\rceil \left(\left\lceil \frac{k}{2} \right\rceil+1 \right)\\
&=  \left( \frac{1}{2}\left\lceil \frac{k}{2} \right\rceil -\frac{1}{2}\right) \left(s \left\lceil \frac{k}{2} \right\rceil + r\left(\left\lceil \frac{k}{2} \right\rceil +1\right)\right)+\frac{r}{2}\left(\left\lceil \frac{k}{2} \right\rceil+1\right)\\
& \ge \left( \frac{1}{2}\left\lceil \frac{k}{2} \right\rceil -\frac{1}{2}\right)  (n-(t-1)(\lceil {k}/{2}\rceil-1)). 
\end{align*} This proves Theorem~\ref{e(H)}(\ref{eD}).\medskip

To prove Theorem~\ref{e(H)}(\ref{conH}),    suppose that $H$ is disconnected. Let $x, y \in V(H)$ be such that $x$ and $y$ are in different components of $H$. By Theorem~\ref{e(H)}(\ref{NHu}),   $\{x, y\} \subseteq D_i$ for some $i \in [p]$, and there must exist a vertex $w \in D_j$ such that $xw \not\in E(H)$ and $yw \in E(H)$, where   $j \in [p]$ with $j \ne i$. By Theorem~\ref{e(H)}(\ref{NHu}),  $x$ and $w$ have at least $t-2$ common neighbors in $H$. But then $x$ and $y$ must be in the same component of $H$, a contradiction. This proves Theorem~\ref{e(H)}(\ref{conH}). \medskip

It remains to  prove Theorem~\ref{e(H)}(\ref{eqH}). By Theorem~\ref{e(H)}(\ref{conH}), $H$ is connected. Let   $q\in\mathbb N$ with $q \ge t-1$. Assume $\delta(H) \ge q$. Following Day    \cite{Day2017}, we next apply the $q$-neighbour    bootstrap percolation on $H$. Note that  $H$ is not necessarily $K_t$-saturated. Given a set  $S \subseteq V(H)$ and  any vertex  $v\in V(H)$, let  $N_S(v): =N_H(v) \cap S$ and $d_S(v): =|N_S(v)|$. 
Let  $R \subseteq V(H)$ be  any nonempty set.   Let $R^0: =R$ and for   $i \ge 1$, let 
\[
R^i: =R^{i-1} \cup \{v \in V(H):  d_{R^{i-1}}(v) \ge q\}.
\] 
 Let $\overline{R}: = \bigcup_{i \ge 0} R^i$, the closure of $R$ under the $q$-neighbor bootstrap percolation on $H$.   Then 
 \[e(H[\overline{R}]) \ge q(|\overline{R}|-|R|), \] 
 because every vertex in $R^i\less R^{i-1}$ is 
 adjacent to at least $q$ vertices in $R^{i-1}$.  
 Let  $Y(R): =V(H) \less \overline{R}$. Finally, for any $v \in V(H)$,  let 
  \[
  \omega_{_R}(v): =d_{\overline{R}}(v) + d_{Y(R)}(v)/2.
  \]
  We call $\omega_{_R}(v)$ the  weight of $v$ (with respect to $R$). Then 
  \[e_H(\overline{R}, Y(R)) +e(H[Y(R)])= \sum_{v \in Y(R)} \omega_{_R}(v).\]
Within $Y(R)$, we define $B(R)$ to be the set $\{v \in Y(R):  \omega_{_R}(v) < q\}$, which we call the set of bad vertices.   
We next show that  there exists a constant $c_1(q, k)$  and a nonempty set $R\subseteq V(H)$ with $|R|\le c_1(q, k)$ such that $B(R)=\es$.   \medskip 
   
  Assume $B(R)\ne \es$ for our initial $R$. Our goal is  to  move a small number of vertices into $R$ so that the remaining vertices in
$B(R)$ have strictly larger weight.  To achieve this,   let 
\[\mathcal{U}_R: =\{U \subseteq R:  U=N_R(v)\ \text{for some}\ v \in B(R)\}.\]
 Note that for every vertex  $v\in B(R)$, $d_R(v)\le q-1$.  Thus 
 \[|\mathcal{U}_R|\le 1+|R|+{|R| \choose 2}+ \cdots+{|R| \choose q-1}.\]
   Let  $\mathcal{U}_R:=\{U_1, U_2, \ldots, U_{|\mathcal{U}_R|}\}$ and let 
  $u_i \in B(R)$ with $N_R(u_i)=U_i$ for all $i \in \{1, 2, \ldots, |\mathcal{U}_R|\}$. Then $d_{\overline{R}}(u_i) <q$, and so $d_{Y(R)}(u_i) \ge 1$ because $d_H(u_i) \ge q$. Let $x_i \in Y(R)$ such that $u_ix_i \in E(H)$ for  all $i \in \{1,  \ldots, |\mathcal{U}_R|\}$, and let $X(R): =\{x_1, x_2, \ldots, x_{_{|\mathcal{U}_R|}}\}$. By the choice of $\mathcal{U}_R$ and $u_1, u_2,\ldots, u_{_{|\mathcal{U}_R|}}$, for every  vertex $v\in B(R)$, we see that $N_R(v)=N_R(u_i)$ for some    $ i \in \{1, 2, \ldots, |\mathcal{U}_R|\}$.  Finally, let 
  \[S(R):=\{v \in B(R):  N_R(v)=N_R(u_i)  \ \text{and}\  \{v, x_i\} \subseteq D_j  \ \text{for some}\ i \in \{1, 2, \ldots, |\mathcal{U}_R|\} \ \text{and }\ j \in [p] \}.\]
   We next show that  {\bf Algorithm}~\ref{algo} below yields   a nonempty set   $R\subseteq V(H)$ with $B(R)=\emptyset$.\bigskip

\begin{algorithm}[H]\label{algo}
\SetAlgoLined
\KwData{  $H:=G \less (\bigcup_{i \in [p]} E(G[V(D_i)]))$   is a spanning subgraph of $G_r$ with  $\delta(H) \ge q$  }
\KwResult{A nonempty set $R\subseteq V(H)$ with $B(R)=\emptyset$}
Set  $R $  to be a set containing an arbitrary  vertex in $H$;

\While{$B(R)\neq \emptyset$ }{
Set $R$ to be $R \cup X(R) \cup S(R) \cup \bigcup_{j=1}^{|\mathcal{U}_R|}N_{\overline{R}}(x_j)$\;
}

\caption{Building a nonempty set $R\subseteq V(H)$ with $B(R)=\emptyset$}
\end{algorithm} 
\bigskip

\noindent Let $R_i$ be the set $R$ obtained  in the $i$-th iteration  of {\bf Line 2}  when  running {\bf Algorithm}~\ref{algo}. Then for all $i\ge1$, $R_{i-1}\subseteq R_i$,    $\overline{R}_{i-1} \subseteq \overline{R}_i$, $Y(R_i)\subseteq Y(R_{i-1})$ and $B(R_i) \subseteq B(R_{i-1})$. To see why $\omega_{_{R_i}}(v)>\omega_{_{R_{i-1}}}(v)$ for all $v\in  B(R_i)$, we next introduce a control function on $V(H)$, because   dealing with $\omega_{_R}(v)$ directly is difficult.   Let   $\phi_{_R}(v): =\sum_{x \in N_H(v)} f_{_R}(x)$ for all $v\in V(H)$, where for all $x\in V(H)$, 
$$f_{_R}(x)=
\begin{cases}
1, & \text{if}\ x \in R, \\
1/2, & \text{if}\ x \in \overline{R} \less R,\\
d_{R}(x)/(2q), & \text{if}\ x \in Y(R). 
\end{cases}
$$

It is worth noting that  $\phi_{_R}(v) \le \omega_{_R}(v)$ for every vertex  $v \in V(H)$, because  $d_{Y(R)}(x)\ge1$ and  $d_R(x) \le q-1$ for all  $x \in Y(R)$. Similarly,  for all $i\ge1$, $f_{_{R_{i-1}}}(x) \le f_{_{R_i}}(x)$ for every $x \in V(H)$,  because $Y(R_i)\subseteq Y(R_{i-1})$.  We next claim that \bigskip

\noindent ($\ast$)  for all $i \ge 1$ and  every $v \in B(R_i)$,  $\phi_{_{R_i}}(v) \ge \phi_{_{R_{i-1}}}(v)+1/(2q)$.

\pf Let $i\ge1$ and $v \in B(R_i)$. Then $v\in B(R_{i-1})$, since $B(R_i) \subseteq B(R_{i-1})$. Let  $\mathcal{U}_{R_{i-1}}$, $\{u_1,   \ldots, u_{_{|\mathcal{U}_{R_{i-1}}|}}\}\subseteq B(R_{i-1})$, and  $\{x_1,   \ldots, x_{_{|\mathcal{U}_{R_{i-1}}|}}\}\subseteq Y(R_{i-1})$  be defined accordingly for $R_{i-1}$.  Then  $N_{R_{i-1}}(v)=N_{R_{i-1}}(u_j)$ for some $j \in \{1,2, \ldots, |\mathcal{U}_{R_{i-1}}|\}$. To prove $\phi_{_{R_i}}(v) \ge \phi_{_{R_{i-1}}}(v)+1/(2q)$, it suffices to show that $f_{_{R_i}}(x) \ge f_{_{R_{i-1}}}(x)+1/(2q)$ for some $x \in N_H(v)$. Since $\{x_1,   \ldots, x_{_{|\mathcal{U}_{R_{i-1}}|}}\}\subseteq Y(R_{i-1})\cap R_i$,  we see that $f_{_{R_{i-1}}}(x )=d_{_{R_{i-1}}}(x )/(2q) \le (q-1)/(2q)=1/2-1/(2q)$,  and $f_{_{R_i}}(x )=1>f_{_{R_{i-1}}}(x )+1/(2q)$ for all $x\in \{x_1,   \ldots, x_{_{|\mathcal{U}_{R_{i-1}}|}}\}$.  We may assume that $vx_j \not\in E(H)$ for all $j\in \{1, \ldots, |\mathcal{U}_{R_{i-1}}|\}$, otherwise we are done. Since $v\in B(R_i)$, by the choice of $x_j$ and  $S(R_{i-1})$, we see that $\{v, x_j\} \nsubseteq V(D_\ell)$ for all $\ell \in [p]$. By Theorem~\ref{e(H)}(\ref{NHu}) applied to $v$ and $x_j$, $H[N_H(v) \cap N_H(x_j)]$ has a $K_{t-2}$ subgraph. Let $W$ be the vertex set of such a $K_{t-2}$ subgraph.  It follows that $W \nsubseteq R_{i-1}$, else $G_r[W \cup \{u_j, x_j\}]=K_t$, since $N_{R_{i-1}}(v)=N_{R_{i-1}}(u_j)$ and  $u_jx_j \in E(H)$. Let $x \in W \less R_{i-1}$.   \medskip

  If $x \in \overline{R}_{i-1} \less R_{i-1}$, then $f_{_{R_{i-1}}}(x)=1/2$ and $f_{_{R_i}}(x)=1$, and so $f_{_{R_i}}(x) \ge f_{_{R_{i-1}}}(x)+1/(2q)$, as desired.   If $x \in Y(R_{i-1})$, then either $x \in \overline{R}_i$ or $x \in Y(R_i)$. In both cases, we have  $f_{_{R_{i-1}}}(x)=d_{R_{i-1}}(x)/(2q)\le 1/2-1/(2q)$.  If  $x \in \overline{R}_i$, then $f_{_{R_i}}(x) \ge 1/2$ and so $f_{_{R_i}}(x) \ge f_{_{R_{i-1}}}(x)+1/(2q)$.   Finally, if $x \in Y(R_i)$, then $d_{R_i}(x)\ge d_{R_{i-1}}(x)+1$ because $x_j\in R_i\less R_{i-1}$ and $R_{i-1}\subseteq R_i$. Hence,  $f_{_{R_i}}(x)=d_{_{R_i}}(x)/(2q) \ge (d_{_{R_{i-1}}}(x)+1)/(2q)= f_{_{R_{i-1}}}(x)+1/(2q)$.  \medskip 
  
  In all cases, we have shown that there exists  some  vertex $x\in N_H(v)$ such that     $f_{_{R_i}}(x)\ge f_{_{R_{i-1}}}(x)+1/(2q)$.  Hence,  $\phi_{_{R_i}}(v) \ge \phi_{_{R_{i-1}}}(v)+1/(2q)$ for all $i\ge1$ and $v\in B(R_i)$.   \hfill\vrule height3pt width6pt depth2pt\\

\noindent By ($\ast$),    {\bf Algorithm}~\ref{algo} stops  after   $m\le 2q^2$ iterations of {\bf Line 2}. Hence  $R_m\subseteq V(H)$ with $  R_m\ne\es$ but  $B(R_m)=\emptyset$. For all $i \ge 0$, 
\begin{align*}
|R_{i+1}| &= |R_i|+|X(R_i)| +|S(R_i)| + |\bigcup_{_{j=1}}^{|\mathcal{U}_{R_i}|}N_{\overline{R}_i}(x_j)|\\
& \le |R_i|+|\mathcal{U}_{R_i}|+(k-2)|\mathcal{U}_{R_i}|+(q-1)|\mathcal{U}_{R_i}|\\
&=|R_i|+(k+q-2)|\mathcal{U}_{R_i}|\\
&\le |R_i|+ (k+q-2)\left(1+|R_i|+{|R_i| \choose 2}+\cdots+{|R_i| \choose q-1}\right), 
\end{align*}
 which   only depends on $q$ and $k$. It follows that by {\bf Algorithm}~\ref{algo},  there exists  a constant $c_1(q,k)$  and a non-empty set $R\subseteq V(H)$ with   $|R| \le c_1(q,k)$  such that  $B(R)=\emptyset$.  Then $\omega_{_R}(v)\ge q$ for all $v \in Y(R)$ and so  
  \[e_H(\overline{R}, Y(R)) +e(H[Y(R)])= \sum_{v \in Y(R)} \omega_{_R}(v)\ge q|Y(R)|.\]
Therefore, 
\begin{align*}
  e(H) & = e(H[\overline{R}])+e_H(\overline{R}, Y(R)) +e(H[Y(R)])\\
  & \ge q(|\overline{R}|-|R|)+q|Y(R)|\\
& \ge q(|\overline{R}|-c_1(q,k))+q|Y(R)|\\
& = q(n-c_1(q,k))\\
&=qn-c(q, k)
\end{align*}
where $c(q, k)=qc_1(q,k)$. This proves Theorem~\ref{e(H)}(\ref{eqH}).\medskip

This completes the proof of Theorem \ref{e(H)}. \hfill\vrule height3pt width6pt depth2pt

\section {Lower bound on the size of $(K_t, \mathcal{T}_k)$-co-critical graphs}\label{Lower} 
  We begin this section with a useful lemma, which may be of independent interest.  It is worth noting that Lemma~\ref{ind} is stronger than Theorem~\ref{Hajnal1} when $\alpha(G)>|G|/2$. We include a proof here for completeness and the proof of Lemma~\ref{ind} is due to Hehui Wu\footnote{The first author   would like to  thank Hehui Wu,  from Fudan University in  China,   for helpful discussion  on the statement of Lemma~\ref{ind} and his neat proof of Lemma~\ref{ind}, which   is completely different from the one of Hajnal  \cite{Hajnal}.}.  For a graph $G$, a set $A\subseteq V(G)$ is \dfn{stable} if $G[A]$ has no edges. We use $\alpha(G)$ and $\omega(G)$ to denote the independence number and clique number of   $G$, respectively. 
   
 \begin{thm}[Hajnal  \cite{Hajnal}]\label{Hajnal1}
Let $G$ be a graph and let $\mathcal{F}$ be the family of all maximum stable sets of $G$. Then \[
\left| \bigcap_{S \in \mathcal{F}} S\right|+\left| \bigcup_{S \in \mathcal{F}} S\right| \ge  2 \alpha(G).\] 
\end{thm}
 
\begin{lem}\label{ind}
Let $G$ be a graph with $\alpha(G)>|G|/2$ and let $\mathcal{F}$ be the family of all maximum stable sets of $G$. Then 
\[\left| \bigcap_{S \in \mathcal{F}} S\right| \ge  \delta(G)+2\alpha(G)-|G|\ge\delta(G)+1.\] Moreover,  if $ \bigcap_{S \in \mathcal{F}} S=\{u\}$, then  $\alpha(G)=(|G|+1)/2$ and $u$ is an isolated vertex in $G$. 
\end{lem}

\pf  Let $X \in \mathcal{F}$ and $Y:=V(G) \less X$. Then $|X|=\alpha(G)>|G|/2$, and so $|X| > |Y|$. Let $H:=G[X, Y]$ be the  bipartite subgraph of $G$ with   $V(H)=X\cup Y$ and   $E(H)=\{xy\in E(G): x\in X, y\in Y\}$. Let $T$ be a maximum stable set of $H$ and  let $X_1:=X \less T$,  $Y_1:=Y \cap T$ and $Y_2:=Y \less T$.  Then $|Y_1|+|X\less X_1|=|T| \ge |X| =|X_1|+|X\less X_1|> |Y|=|Y_1|+|Y_2|$,   which implies that  $|X_1| \le |Y_1|$ and $|X\less X_1| >|Y_2| $. We next show that $H':=G[X\less X_1, Y_2]$ contains a matching that saturates $Y_2$. For any $S \subseteq Y_2$, we have $|N_{H'}(S)| \ge |S|$, else $T':=(T \less N_{H'}(S)) \cup S$ is a stable set of $H$ with $|T'| >|T|$, a contradiction. By Hall's Theorem, there exists a matching, say  $M$,  of $H'$ that saturates $Y_2$.  Let $X_2 := V(M) \cap X$ and $X_3:=X \less (X_1 \cup X_2)$. Then \[
|X_3|=|X|-|X_1|-|X_2| \ge |X|-|Y|=2\alpha(G)-|G|>0,\]
 because $|X_1| \le |Y_1|$,  $|X_2|=|Y_2|$ and $\alpha(G)>|G|/2$. 
Note that $X_1 \cup Y_1$ is anti-complete to $X \less X_1$ in $H$.  By the choice of $T$,   $\alpha(H[X_1 \cup Y_1]) \le |X_1|$. Moreover, $\alpha(H[X_2 \cup Y_2]) \le |X_2|$ because $M$ is a perfect matching  of $G[X_2,  Y_2]$. Then for any $S \in \mathcal{F}$, $|S \cap (X_1 \cup Y_1)| \le |X_1|$ and $|S\cap (X_2 \cup Y_2)| \le |X_2|$. Therefore, $|X_3| \ge |S \cap X_3| =|S| - |S \cap (X_1 \cup Y_1)|- |S \cap (X_2 \cup Y_2)| \ge |X| - |X_1|-|X_2|=|X_3|$. It follows that $|S \cap X_3| = |X_3|$.  Then $X_3 \subseteq  S$. Hence,   $X_3 \subseteq \bigcap_{S \in \mathcal{F}} S$ by the arbitrary choice of $S$.  \medskip

Next, suppose   there exists a vertex  $u \in X_3$ with $d_G(u)=d>0$. Let   $N_G(u):=\{ v_1, \ldots, v_d \}$. Then $\{ v_1, \ldots, v_d \} \subseteq Y_2$.  Let $u_1, \ldots, u_d \in X_2$ be such that $u_iv_i \in E(M)$ for all $i \in [d]$.  For each $i \in [d]$, let $M^i:=(M\less u_iv_i)\cup  \{uv_i\}$, $X_2^i:=V(M^i) \cap X$ and $X_3^i:= X \less (X_1 \cup X_2^i)$. Then $u_i \in X_3^i$ and $M^i$ is a perfect matching  of $G[X^i_2,  Y_2]$.  By the arbitrary choice of $M$,   $u_i\in  \bigcap_{S \in \mathcal{F}} S$. Therefore,  $|\bigcap_{S \in \mathcal{F}} S| \ge | \{u_1, \ldots, u_d\}\cup X_3  | \ge d+(2\alpha(G)-|G|) \ge \delta(G)+2\alpha(G)-|G|\ge\delta(G)+1$, as desired. \medskip

Finally, if   $\bigcap_{S \in \mathcal{F}} S = \{u\}$, then $1=|\bigcap_{S \in \mathcal{F}} S| \ge d+2\alpha(G)-|G|$.   It follows that $d=0$ and $\alpha(G)=(|G|+1)/2$, because $2\alpha(G)-|G| >0$. This completes the proof of Lemma~\ref{ind}. 
\hfill\vrule height3pt width6pt depth2pt \\

 We are now ready to prove Theorem~\ref{lower}. \\

\noindent {\bf Proof of Theorem~\ref{lower}:} Let $G$ be a  $(K_t, \mathcal{T}_k)$-co-critical graph on  $n $   vertices, where  $t \ge 4$ and  $k \ge \max\{6,t\}$.  Then $n \ge (t-1)(k-1)+1$ and     $G$ admits  a critical-coloring. Among all critical-colorings of $G$, let $\tau: E(G) \rightarrow \{$red, blue$\}$ be a critical-coloring  of $G$ with $|E_r|$   maximum. By the choice of $\tau$, $G_r$ is $K_t$-saturated and $G_b$ is  $  \mathcal{T}_k$-free. By Theorem~\ref{e(H)}(a),    $\delta(G_r)\ge 2t-4$.   
  Let $D_1, \cdots, D_p$ be all components of $G_b$. By Lemma~\ref{structure}(a),  $|D_i|\le k-1$ for all $i\in [p]$. Then $ (t-1)(k-1)+1\le n\le p(k-1)$. This implies that   $p \ge t$. Let $H:=G \less (\bigcup_{i \in [p]} E(G[V(D_i)]))$. Then $H$ is a spanning subgraph of $G_r$. Clearly, $H$ is $K_t$-free.  \medskip

  Assume first that $\delta(H) \ge 2t-4$. By  Theorem~\ref{e(H)}(\ref{eqH}) applied to $H$ and $q=2t-4$,  there exists a constant $c(2t-4,k)$ such that $e(H) \ge (2t-4)n-c(2t-4,k)$. This, together with     Theorem~\ref{e(H)}(\ref{eD}), yields that 
\begin{align*}
e(G)&=e(H)+\sum_{i=1}^p e(G[V(D_i)]) \\
&\ge (2t-4)n-c(2t-4,k) + \left( \frac{1}{2}\left\lceil \frac{k}{2} \right\rceil -\frac{1}{2}\right)(n-(t-1)(\lceil k/2 \rceil-1)) \\
&= \left(\frac{4t-9}{2}+\frac{1}{2}\left\lceil \frac{k}{2} \right\rceil \right)n-c(2t-4,k)- \frac{1}{2}(t-1)(\lceil k/2 \rceil-1)^2\\
&=\left(\frac{4t-9}{2}+\frac{1}{2}\left\lceil \frac{k}{2} \right\rceil \right)n-c_1(t,k),
\end{align*}
as desired, where $c_1(t, k)= c(2t-4,k)+ \frac{1}{2}(t-1)(\lceil k/2 \rceil-1)^2$. 
\medskip

Assume next that $\delta(H) \le 2t-5$.  Note that  $k \ge \max\{6,t\}\ge t$ for all $t \ge 4$.  Let  $u\in V(H)$ with $d_H(u)=\delta(H) $. We may assume that $u \in V(D_p)$. Let $N_H(u)=\{u_1, \ldots,  u_{\delta(H)}\}$. By Theorem~\ref{e(H)}(\ref{NHu}) applied to $u$ and any vertex in  $V(H) \less (V(D_p) \cup N_H(u))$,   we see that $ H[N_H(u)] $ must have  a $K_{t-2}$ subgraph.  We may  assume that $  H[\{u_1, \ldots, u_{t-2}\}]=K_{t-2} $. Then we may further assume that $u_i \in V(D_{p-i})$ for all $i \in [t-2]$. Note that $ H[N_H(u)] $ is $K_{t-1}$-free and $\omega( H[N_H(u)])=t-2>| N_H(u)|/2$. Let $\mathcal{F}$ be the family of all $K_{t-2}$ subgraphs of $ H[N_H(u)] $. 
By Theorem~\ref{e(H)}(\ref{nocom}), $ |\bigcap_{A\in \mathcal{F}}A|\le1$. 
  By Lemma~\ref{ind} applied to the complement  of $ H[N_H(u)]$, we have $ |\bigcap_{A\in \mathcal{F}}A|=1$. We may assume that $\bigcap_{A\in \mathcal{F}}A=\{u_1\}$.  By Lemma~\ref{ind} again, $|N_H(u)|=2t-5$,     $u_1$ is complete to $N_H(u) \less u_1$ in $H$ and  $u_1$ is contained in all $K_{t-2}$ subgraphs of $H[N_H(u)]$. Then  $H[N_H(u)]\less u_1$ is $K_{t-2}$-free.  By Theorem~\ref{e(H)}(\ref{NHu}) applied to $u$ and any vertex in $V(H) \less (V(D_p) \cup N_H(u))$, $V(H) \less (V(D_p) \cup N_H(u))$ must be complete to  $u_1$ in $H$. Thus $\{u_1\}=V(D_{p-1})$.  By Theorem \ref{e(H)}(\ref{eqH}) applied to $H$ and $q=2t-5$,  there exists a constant $c(2t-5,k)$ such that $e(H) \ge (2t-5)n-c(2t-5,k)$.   Since  $n\ge (t-1)(k-1)+1$ and $|V(D_i)|\le k-1$ for all $i\in[p]$ with $i\ne p-1$, we see that   $p\ge t$.    If  $p=t$, then  $n= (t-1)(k-1)+1$ and $|V(D_i)|= k-1$ for $i\in[p]$ with $i\ne p-1$.  In this case,
  \begin{align*}
e(G)&= e(H)+\sum_{i=1}^{p} e(G[V(D_i)]) \\
&\ge ((2t-5)n-c(2t-5,k)) +  (p-1)(k-1)(k-2)/2\\
&= ((2t-5)n-c(2t-5,k)) +  (n-1)(k-2)/2 \\
&=(2t-6+k/2)n-c(2t-5,k)- (k-2)/2\\
&\ge \left(\frac{4t-9}{2}+\frac{1}{2}\left\lceil \frac{k}{2} \right\rceil \right)n- c_2(t,k) 
\end{align*}
for all   $k \ge 6$, as desired, where $c_2(t,k)=c(2t-5,k)+ (k-2)/2$. \medskip

  Next assume $p\ge t+1$.   Since $k\ge t$, $|N_H(u)|\le 2t-5$, and $G_r$ is $K_t$-free,  by Lemma~\ref{structure}(b), there are   at most $  t-1$ many $D_i$'s satisfying  $u\notin V(D_i)$ and $D_i\less N_H(u)=\es$. We may assume that for all $i\in[p-t]$, $D_1, \ldots, D_{p-t}$ are such that $u\notin V(D_i)$ and $D_i\less N_H(u)\ne\es$.    By Theorem~\ref{e(H)}(\ref{Di}), $|D_i| =k-1$ for all $i\in[p-t]$. Thus  \[\sum_{i=1}^{p} e(G[V(D_i)])\ge (p-t)(k-1)(k-2)/2.\]
  Note that  $n\le (p-1)(k-1)+1$ because $\{u_1\}=V(D_{p-1})$ and $|D_i|\le k-1$ for all $i\in[p]$ with $i\ne p-1$. Therefore,  
 \begin{align*}
e(G)&= e(H)+\sum_{i=1}^{p} e(G[V(D_i)]) \\
&\ge ((2t-5)n-c(2t-5,k)) +  (p-t)(k-1)(k-2)/2\\
&\ge ((2t-5)n-c(2t-5,k)) + \frac{1}{2}\left(\frac{n-1}{k-1}- t+1\right)(k-1)(k-2) \\
&=(2t-6+k/2)n-c(2t-5,k)- (k-2)(tk-t-k+2)/2\\
&\ge \left(\frac{4t-9}{2}+\frac{1}{2}\left\lceil \frac{k}{2} \right\rceil \right)n-c(2t-5,k)- [(t-1)k^2-(3t-4)k+2t-4]/2 \\
&= \left(\frac{4t-9}{2}+\frac{1}{2}\left\lceil \frac{k}{2} \right\rceil \right)n-c_3(t,k) 
\end{align*}
for all   $k \ge 6$, as desired, where $c_3(t,k)=c(2t-5,k)+ [(t-1)k^2-(3t-4)k+2t-4]/2$. \medskip

 Let $\ell(t, k):=\max\{c_1(t,k), c_2(t,k), c_3(t, k)\}$.  This completes the proof of Theorem \ref{lower}.
\hfill\vrule height3pt width6pt depth2pt


\section{Proof of Theorem~\ref{upper}}\label{Upper}
Let $t \in \{4,5\}$,  $k\ge 3$ and $n \ge (2t-3)(k-1)+{\lceil k/2 \rceil}{\lceil k/2 \rceil}-1$. We will construct a  $(K_t, \mathcal{T}_k)$-co-critical graph on $n$ vertices which yields the desired upper bound in Theorem \ref{upper}.\medskip

\begin{figure}[htb]
\centering
\includegraphics[scale=0.8]{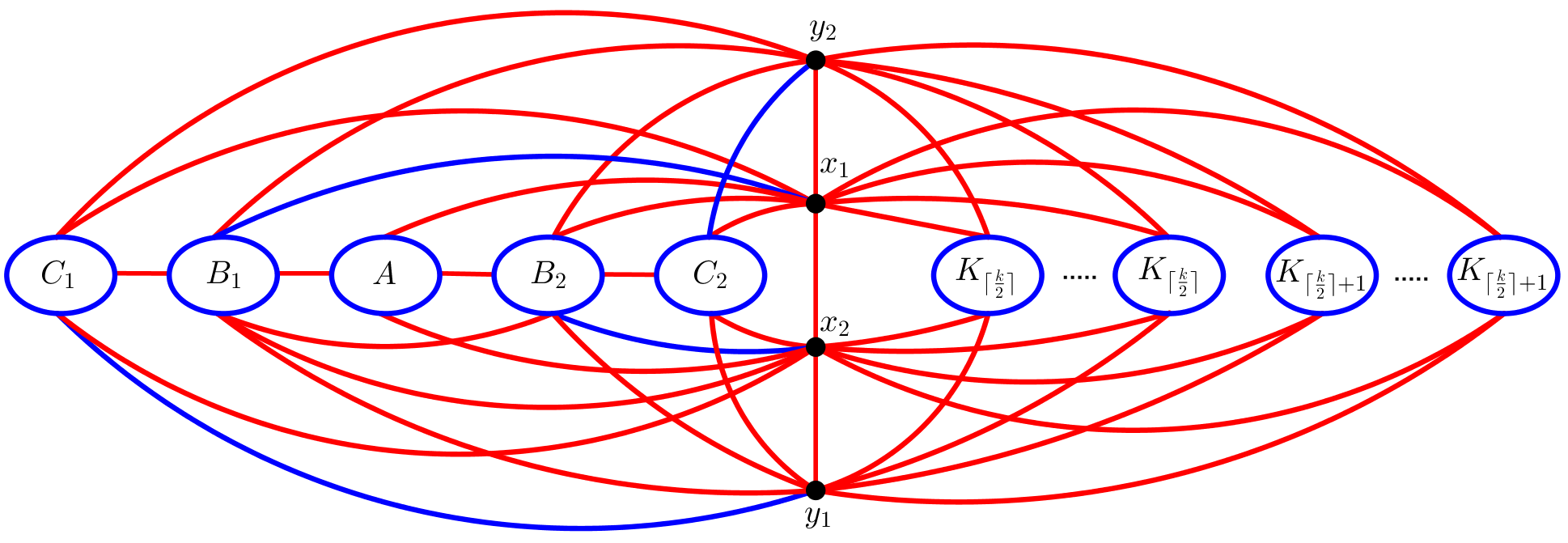}
\caption{A  $(K_4, \mathcal{T}_k)$-co-critical  graph for all $k \ge 4$.}
\label{K4Tk-Saturated}
\end{figure}

\begin{figure}[htb]
\centering
\includegraphics[scale=0.8]{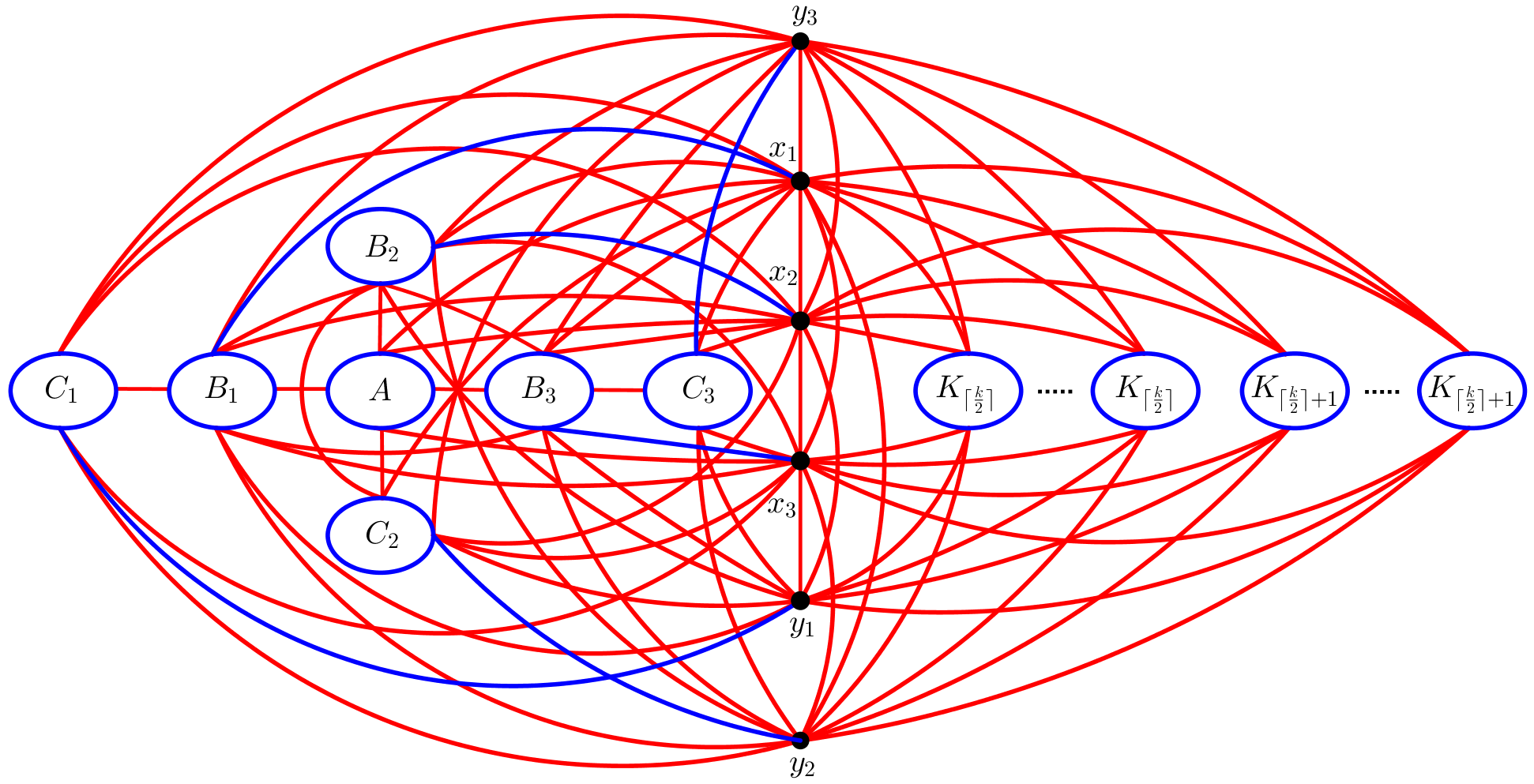}
\caption{A  $(K_5, \mathcal{T}_k)$-co-critical graph for all $k \ge 4$.}
\label{K5Tk-Saturated}
\end{figure}

Let $r, s$ be the remainder and quotient of $n- (2t-3)(k-1)$ when divided by $\lceil k/2 \rceil$, and  let $A: =K_{k-1}$. For each $i \in [t-2]$,  let  $B_i: =K_{k-2}$ and $C_i:=K_{k-2}$.   Let $H_1$ be obtained from disjoint copies of $A, B_1, \ldots, B_{t-2}, C_1, \ldots, C_{t-2}$ by joining every vertex in $B_i$ to all vertices in $A\cup C_i\cup  B_j$ for each   $i\in [t-2]$ and  all $  j \in [t-2]$ with $j \neq i$. Let $H_2:= (s-r)K_{\lceil k/2 \rceil} \cup rK_{\lceil k/2 \rceil+1}$ when $k \ge 4$, and  $H_2: = sK_2 \cup rK_1$ when $k =3$.     Finally, let $G$ be the graph obtained from $H:=H_1 \cup H_2$ by adding $2t-4$ new vertices $x_1, \ldots, x_{t-2}, y_1, \ldots, y_{t-2}$, and then, for each   $i\in [t-2]$, joining: $x_i$ to  every vertex in $V(H)$ and all $x_j$;  and $y_i$ to  every vertex in $V(H) \setminus V(A)$  and all  $x_j$, where  $  j \in [t-2]$ with $j \neq i$. The construction of  $G$ when $t=4$ and $k \ge 4$ is depicted in Figure~\ref{K4Tk-Saturated}, and the construction of  $G$ when $t=5$ and $k \ge 4$ is depicted in Figure~\ref{K5Tk-Saturated}.\medskip


Let  $\sigma: E(G) \rightarrow \{$red, blue$\}$ be defined as follows: all edges in   $A, B_1, \ldots, B_{t-2}, C_1, \ldots, C_{t-2}$ and $H_2$ are colored blue; for every $i\in[t-2]$,  all edges between $x_i$ and $B_i$ are colored blue and  all edges between $y_i$ and $C_i$ are colored blue; the remaining edges of $G$ are all  colored red. Note that the  $\{$red, blue$\}$-coloring of $G$  depicted in Figure \ref{K4Tk-Saturated} (resp.  Figure \ref{K5Tk-Saturated}) is $\sigma$  when $t=4$ (resp. $t=5$) and $k \ge 4$. Clearly,   $\sigma $   is a critical-coloring  of $G$. We next show that $\sigma$ is the unique critical-coloring  of $G$  up to symmetry.  \medskip

Let $X:=\{x_1,\ldots,x_{t-2}\}$ and $Y:= \{y_1, \ldots, y_{t-2}\}$.   Let $\tau: E(G) \rightarrow \{$red, blue$\}$ be an arbitrary critical-coloring of $G$.  It suffices to show that $\tau=\sigma$ upon  to symmetry.  Let $G^\tau_r$ and $G^\tau_b$ be  $G_r$ and $G_b$ under the coloring $\tau$, respectively.  Note that $G[V(A) \cup V(B_1) \cup \cdots \cup V(B_{t-2}) \cup X] = K_{(t-1)(k-1)}$. By Lemma \ref{structure}(a) and the fact that $G^\tau_r$ is $K_t$-free,   $G^\tau_b[V(A) \cup V(B_1) \cup \cdots \cup V(B_{t-2}) \cup X]$ has exactly $t-1$ components, say $D_1, \ldots, D_{t-1}$, such that   $V(D_i)$ is complete to $V(D_j)$ in $G^\tau_r$ for all $i, j \in [t-1]$ with $i \neq j$. Then each $D_i$ is isomorphic to $K_{k-1}$ in $G^\tau_b$ for all $i  \in [t-1]$. Since every vertex in $V(A) \cup V(B_1) \cup \cdots \cup V(B_{t-2}) \cup X$ belongs to a blue $K_{k-1}$ in $G^\tau_b$, it follows  that: for each $i  \in [t-2]$,  $y_i$ is complete to $ V(B_1) \cup \cdots \cup V(B_{t-2}) \cup (X\less x_i)$  in $G^\tau_r$;  and  $V(C_i)$ is complete to $V(B_i)\cup X$ in $G^\tau_r$. We next prove three   claims. \bigskip

\setcounter{counter}{0}
\noindent {\bf Claim\refstepcounter{counter}\label{A}  \arabic{counter}.}  
  $ A= D_i $ for some $i\in [t-1]$. 

\pf Suppose $ A  \neq  D_i $ for all $i\in [t-1]$. Then  for each $i\in [t-1]$,  $(V(B_1) \cup \cdots \cup V(B_{t-2})) \cap V(D_i) \neq \emptyset$. Let $d_i \in (V(B_1) \cup \cdots \cup V(B_{t-2})) \cap V(D_i)$ for all $i\in [t-1]$. Then $d_1, \ldots, d_{t-1}$ are pairwise distinct, and $G^\tau_r[\{d_1, \ldots, d_{t-1} \}]=K_{t-1}$.  But then $G^\tau_r[\{d_1, \ldots, d_{t-1}, y_1\}]=K_t$, because  $y_1$ is complete to $V(B_1) \cup \cdots \cup V(B_{t-2})$ in $G^\tau_r$, a contradiction. This proves that  $ A = D_i $ for some $i\in [t-1]$. \hfill\vrule height3pt width6pt depth2pt\bigskip

 By Claim~\ref{A}, we may assume that $ A = D_{t-1} $. Then $V(A)$ is complete to $ V(B_1) \cup \cdots \cup V(B_{t-2}) \cup X$ in $G^\tau_r$.  For each $i \in [t-2]$, since $G^\tau_b$ is  $\mathcal{T}_k$-free,  there   must exist  a vertex $c_i \in V(C_i)$ such that $c_i$ is adjacent to at most one vertex of $Y$ in $G^\tau_b$.  Then $c_i$ is adjacent to at least $t-3$  vertices of $Y$ in $G^\tau_r$.  We next show that \medskip   
 
\noindent {\bf Claim\refstepcounter{counter}\label{BD1}  \arabic{counter}.} 
   For each  $i \in [t-2]$, $|X \cap V(D_i)|= 1$. 
   
   \pf Suppose $|X \cap V(D_i)|\ne 1$ for some $i\in[t-2]$.  Since $|X|=t-2$, we may assume that  $|X \cap V(D_1)|\ge2$ and $X \cap V(D_{t-2})=\es$.  We may further assume that    $x_1, x_2\in V(D_1)$.  Then $x_1x_2\in E_b$.   Since $X \cap V(D_{t-2})=\es$ and  for all $i\in [t-2]$, $|V(B_i)|=k-2<k-1=|V(D_{t-2})|$, we may assume that $V(B_i)\cap V(D_{t-2})  \ne \es$  for $i\in[2]$.  Let $b_1\in V(B_1)\cap V(D_{t-2})$. We may assume that $c_1y_i\in E_r$ for some $i\in [2]$, because $c_1$ is adjacent to at least $t-3$  vertices of $Y$ in $G^\tau_r$.  If $t=4$, then    $G^\tau_r[\{ b_1,  c_1, y_i, x_{3-i}\}]=K_4$, a contradiction.   Thus $t=5$.   We claim that    $V(B_1)\cap V(D_2)=\es$ and  $V(B_2)\cap V(D_2)=\es$.   Suppose, say  $V(B_1)\cap V(D_2)\ne\es$.  Let $b_2\in V(B_1)\cap V(D_2)$.  Then $G^\tau_r[ \{b_1, b_2,   c_1, y_i, x_{3-i}  \}]=K_5$,  a contradiction. Thus $V(B_1)\cap V(D_2)=\es$ and  $V(B_2)\cap V(D_2)=\es$.   Then $V(D_2)=V(B_3)\cup\{x_3\}$. But  then $G^\tau_r[ \{b_1,   c_1, y_i, x_{3-i}, x_3\}]=K_5$, a contradiction.  \hfill\vrule height3pt width6pt depth2pt\bigskip    
 
\noindent {\bf Claim\refstepcounter{counter}\label{BD}  \arabic{counter}.} 
   For each  $i \in [t-2]$, $V(B_i) \subseteq V(D_j)$ for some  $j \in [t-2]$. 

\pf Suppose there exists an $i\in [t-2]$ such that $V(B_i) \nsubseteq V(D_j)$ for every  $j \in [t-2]$. We may assume $i=1$. Since $V(B_1) \subseteq V(D_1) \cup \cdots \cup V(D_{t-2})$, we see that  $k-2=|B_1| \ge2$. Thus  $k\ge4$.  We claim that $V(B_1)\cap V(D_j) =\es$ for some $j\in [t-2]$. Suppose 
$V(B_1)\cap V(D_j)\ne\es$ for all $j\in [t-2]$.
 Let   $d_j \in V(B_1) \cap V(D_j)$ for all $j\in [t-2]$.  But then $G^\tau_r[\{d_1,  \ldots, d_{t-2}, c_1, y_{\ell}\}] =K_t$, where    $c_1y_\ell\in E_r$ for some $\ell\in[t-2]$,   a contradiction.  Thus $V(B_1)\cap V(D_j) =\es$ for some $j\in [t-2]$, as claimed. We may assume that $V(B_1)\cap V(D_{t-2}) =\es$. Since $V(B_1) \nsubseteq V(D_j)$ for every  $j \in [t-2]$, it follows that  $t=5$,    $V(B_1) \subseteq V(D_1) \cup   V(D_{2})$,  and  $V(B_1)\cap V(D_{1})\ne\es$  and $V(B_1)\cap V(D_{2})\ne\es$.   Let $d_1 \in V(B_1) \cap V(D_1)$ and $d_2 \in V(B_1) \cap V(D_2)$.  By Claim~\ref{BD1}, let $x_i\in X\cap V(D_3)$. Then  $G^\tau_r[\{d_1, d_2, x_i, c_1, y_j\}]=K_5$, where $c_1y_j\in E_r$ for some $j\in [3]$ with $j\ne i$, a contradiction.   
    \hfill\vrule height3pt width6pt depth2pt\bigskip

By  Claim \ref{BD1} and  Claim \ref{BD},    $V(B_i) \cup V(B_j) \nsubseteq D_{\ell}$ for any $i \ne j \in [t-2]$ and all $\ell \in [t-2]$. 
  By symmetry,    we may assume that $V(B_i) \subseteq V(D_i)$ for all $i \in [t-2]$. Then $V(B_i) \cup \{x_j\} = V(D_{i})$ for some $j  \in [t-2]$ since $|V(D_{i})|=|V(B_i)|+1$ and $V(B_1) \cup \cdots \cup V(B_{t-2})\cup X=V(D_1) \cup \cdots \cup V(D_{t-2})$. 
By symmetry, we may assume that $V(B_i) \cup \{x_i\} = V(D_{i})$ for all $i \in [t-2]$. It follows that for all $i, j \in [t-2]$ with $i \neq j$, $B_i$ is complete to $B_j$ in $G^\tau_r$, $x_i$ is complete to $X \less x_i$ and $B_j$ in $G^\tau_r$, $y_i$ is complete to $C_i$ in $G^\tau_b$, $y_i$ is complete to $C_j \cup (X \less x_i)$ in $G^\tau_r$, $x_i$ is complete to $B_i$ in $G^\tau_b$, $\{x_i, y_i\}$ is complete to $H_2$ in $G^\tau_r$, all edges in $A, B_1, \ldots, B_{t-2}, C_1, \ldots, C_{t-2}$ and $H_2$ are colored blue under $\tau$. This proves that $\tau=\sigma$ and thus $\sigma$ is the unique critical-coloring of $G$ upon to symmetry. It can be easily checked that adding any edge $e\in E(\overline{G})$ to $G$ creates a red $K_t$ if $e$ is colored red, and a blue $T_k$ if $e$ is colored blue. Hence, $G$ is $(K_t, \mathcal{T}_k)$-co-critical.  Note that $e_G(X\cup Y, V(G)\less (X\cup Y))=(t-2)(n-(2t-4))+(t-2)(n-(2t-4+k-1))=(t-2)(2n-4t-k+9)$; $e(G[X\cup Y])= {t-2 \choose 2}+(t-2)(t-3)$; $e_G(V(B_1) \cup \cdots \cup V(B_{t-2}), V(C_1) \cup \cdots \cup V(C_{t-2}))=(t-2)(k-2)^2$; $e(G[V(C_1) \cup \cdots \cup V(C_{t-2})])=(t-2) {k-2\choose 2}$; $e(G[V(A)\cup V(B_1) \cup \cdots \cup V(B_{t-2})])={(t-2)(k-2)+k-1 \choose 2}$. 
Using the facts that $s  \lceil k/2 \rceil + r =n-(2t-3)(k-1)$ and $r \le \lceil k/2 \rceil -1$,  we see that
\begin{align*}
e(G)&= (t-2)(2n-4t-k+9)+{t-2 \choose 2}+(t-2)(t-3)+(t-2)(k-2)^2 \\
&\quad +(t-2) {k-2\choose 2}+{(t-2)(k-2)+k-1 \choose 2} +(s-r){\lceil k/2 \rceil\choose 2} + r {\lceil k/2 \rceil +1 \choose 2}\\
&=(2t-4)n-(t-2)k-\frac{1}{2}(t-2)(5t-9)\\
&\quad+(k-2) \big((t-2)(k-2)+(t-2)(k-3)/2+(t-1)(tk-k-2t+3)/2 \big) \\
&\quad +\frac{s-r}{2} \left\lceil \frac{k}{2} \right\rceil  \left(\left\lceil \frac{k}{2} \right\rceil-1 \right)+\frac{r}{2} \left\lceil \frac{k}{2} \right\rceil \left(\left\lceil \frac{k}{2} \right\rceil+1 \right)\\
&=(2t-4)n-(t-2)k-\frac{1}{2}\big(t-2)(5t-9)+\frac{1}{2}(k-2)((t^2+t-5)k-2t^2-2t+11\big) \\
&\quad + \frac{1}{2}\left( \left\lceil \frac{k}{2} \right\rceil -1\right) \left(s \left\lceil \frac{k}{2} \right\rceil + r\right)+\frac{r}{2}\left(\left\lceil \frac{k}{2} \right\rceil+1\right)\\
&\le (2t-4)n+\frac{1}{2}\big((t^2+t-5)k^2-(4t^2+6t-25)k -t^2+23t-40  \big)\\
&\quad + \frac{1}{2}\left( \left\lceil \frac{k}{2} \right\rceil -1 \right)(n-(2t-3)(k-1)) + \frac{1}{2}\left(\left\lceil \frac{k}{2} \right\rceil-1 \right) \left(\left\lceil \frac{k}{2} \right\rceil+1 \right) \\
&= \left(\frac{4t-9}{2}+\frac{1}{2} \left\lceil \frac{k}{2} \right\rceil \right)n
+\frac{1}{2}(t^2+t-5)k^2-(2t^2+2t-11)k \\
&\quad -\frac{(t-2)(t-19)}{2} 
- \frac{1}{2} \left\lceil \frac{k}{2} \right\rceil \left((2t-3)(k-1) - \left\lceil \frac{k}{2} \right\rceil \right)\\ 
&=\left(\frac{4t-9}{2}+\frac{1}{2} \left\lceil \frac{k}{2} \right\rceil \right)n+C(t,k),
\end{align*}
where $C(t,k)=\frac{1}{2}(t^2+t-5)k^2-(2t^2+2t-11)k  -\frac{(t-2)(t-19)}{2} 
- \frac{1}{2} \left\lceil \frac{k}{2} \right\rceil \left((2t-3)(k-1) - \left\lceil \frac{k}{2} \right\rceil  \right)$.  \medskip 
 
 This completes the proof of Theorem~\ref{upper}. \hfill\vrule height3pt width6pt depth2pt

\end{document}